\documentclass{article}
\usepackage{amsmath, amsthm, amssymb, mathtools}
\usepackage[hidelinks]{hyperref}

\theoremstyle{plain}
\newtheorem{theorem}{Theorem}
\newtheorem{lemma}[theorem]{Lemma}
\newtheorem{proposition}[theorem]{Proposition}
\newtheorem{corollary}[theorem]{Corollary}
\theoremstyle{definition}
\theoremstyle{remark}
\newtheorem{remark}[theorem]{Remark}

\title{First-Order Asymptotics of Greenwood's Variance Estimator}

\date{}
\author{J{\"u}rgen Rodenkirchen\\[0.4em]
\normalsize Medizinische Fakult{\"a}t OWL\\
\normalsize Biostatistik und Medizinische Biometrie\\
\normalsize Universit{\"a}t Bielefeld\\
\normalsize Universit{\"a}tsstra{\ss}e 25\\
\normalsize 33615 Bielefeld\\[0.3em]
\normalsize\href{mailto:Juergen.rodenkirchen@uni-bielefeld.de}{\texttt{Juergen.rodenkirchen@uni-bielefeld.de}}}

\hypersetup{
  pdftitle={First-Order Asymptotics of Greenwood's Variance Estimator},
  pdfauthor={J{\"u}rgen Rodenkirchen}
}

\begin{document}
\maketitle

\begin{abstract}
Let \(\widehat S_n(t)\) denote the Kaplan--Meier estimator and write Greenwood's estimator of
its variance as \(\widehat G_n(t)=\widehat S_n(t)^2\widehat H_n(t)\), where \(\widehat H_n(t)\)
denotes the sum in Greenwood's formula. We refer to this term as the Greenwood component. For
fixed \(t\), \(\widehat H_n(t)=O_p(n^{-1})\), so the natural scaled component is \(\widetilde
H_n(t)=n\widehat H_n(t)\).

Its consistency follows from existing large-sample theory. We study its pointwise first-order
fluctuations about \(h(t)=\int_0^t \lambda(s)/y(s)\,ds\), where \(y(s)=P(X\ge s)\), at the
\(\sqrt n\)-scale. Under independent right censoring, positivity of the at-risk probability, a
continuous-time no-ties model, and mild hazard regularity conditions, we derive the
asymptotically linear representation \(\sqrt n\{\widetilde
H_n(t)-h(t)\}=n^{-1/2}\sum_{i=1}^n\phi_{H,t}(Z_i)+o_p(1)\), with influence function
\(\phi_{H,t}\) specified in Theorem 7.

It has two first-order contributions: the martingale fluctuation of the event counting process
and a contribution generated by fluctuations of the empirical at-risk process. Consequently,
\(\sqrt n\{\widetilde H_n(t)-h(t)\}\xRightarrow{d} N\{0,v_H(t)\}\), with
\(v_H(t)=\operatorname{Var}\{\phi_{H,t}(Z)\}\); the risk-set contribution is part of the full
asymptotic variance and cannot in general be omitted.

Combining this representation with the Kaplan--Meier influence function yields a joint central
limit theorem for \((\widehat S_n(t),\widetilde H_n(t))\), retaining their asymptotic
cross-covariance. Applying the multivariate delta method to \(g(S,h)=S^2h\) gives the
first-order asymptotic distribution of the scaled Greenwood estimator \(\widetilde
G_n(t)=n\widehat G_n(t)=\widehat S_n(t)^2\widetilde H_n(t)\).

These results enable corresponding asymptotic variance calculations for Greenwood's variance
estimator.
\end{abstract}

\section{Introduction}

Greenwood's formula is the classical variance estimator for the
Kaplan--Meier product-limit estimator.  For fixed \(t\), it can be written in
the product form
\[
\widehat G_n(t)
=
\widehat S_n(t)^2\,\widehat H_n(t),
\]
Here \(\widehat S_n(t)\) denotes the Kaplan--Meier estimator, and we refer to
\[
\widehat H_n(t)
=
\sum_{t_j\le t}
\frac{d_j}{n_j(n_j-d_j)},
\qquad n_j:=Y_n(t_j-),
\]
as the Greenwood component.  This decomposition is particularly useful when
the sampling variability of Greenwood's estimator itself is of interest: it
separates the familiar product-limit estimator from the component that carries
the additional random-denominator structure of Greenwood's formula.

The first issue is one of scale.  The pointwise variance of the Kaplan--Meier
estimator is of order \(n^{-1}\), suggesting that Greenwood's variance estimator
should naturally be studied on the same scale.  This scaling can also be
read directly from existing product-limit asymptotic theory. In particular,
Gill \cite[Theorem~4.2.2]{gill1980} gives the weak limit of the
product-limit estimator together with a consistent estimator of the variance
function of its Gaussian martingale factor. In the present continuous-time
no-ties notation, that variance-function estimator is precisely
\[
n\widehat H_n(t)=\widetilde H_n(t).
\]
Consequently, \(n\widehat H_n(t)=O_p(1)\), and hence
\(\widehat H_n(t)=O_p(n^{-1})\).  Since
\(0\le\widehat S_n(t)\le1\), the product representation then gives
\[
\widehat G_n(t)
=
\widehat S_n(t)^2\widehat H_n(t)
=
O_p(n^{-1}).
\]
The natural \(O_p(1)\)-scaled quantities are therefore
\[
\widetilde H_n(t):=n\widehat H_n(t),
\qquad
\widetilde G_n(t):=n\widehat G_n(t)
=
\widehat S_n(t)^2\widetilde H_n(t).
\]

Large-sample theory for the Kaplan--Meier estimator is classical. Early weak-convergence
and efficiency results include Breslow and Crowley \cite{breslow1974}, Wellner
\cite{Wellner1982}, and Gill \cite{gill1983}. Functional and asymptotic
representations were developed from several complementary viewpoints, including
Peterson \cite{peterson1977}, Reid \cite{reid1981}, and Lo and Singh
\cite{loSingh1986}; strong uniform consistency was studied by F{\"o}ldes and
Rejt{\H{o}} \cite{foldesRejto1981}, and later central limit theory for
Kaplan--Meier integrals was developed by Stute \cite{stute1995}. Standard
monographic treatments include Kalbfleisch and Prentice \cite{kalbfleisch1980},
Fleming and Harrington \cite{fleming1991}, and Andersen et al.\
\cite{andersen1993}. For the Greenwood component, consistency of the
appropriately scaled quantity is therefore already contained in Gill's
product-limit theory \cite[Theorem~4.2.2]{gill1980}. The question addressed
here begins at the next asymptotic order: what are the \(\sqrt n\)-fluctuations of
\(\widetilde H_n(t)\) around its deterministic limit
\[
h(t)=\int_0^t \frac{\lambda(s)}{y(s)}\,ds,
\]
and how do these fluctuations combine with those of the Kaplan--Meier
estimator in Greenwood's product representation?

This normalization is also natural for the motivating log-scale variance
calculations.  For any positive random quantity \(X_n\), multiplication by the
deterministic factor \(n\) changes its logarithm only by the constant \(\log n\),
so that
\[
\operatorname{Var}\{\log(nX_n)\}
=
\operatorname{Var}\{\log X_n\}.
\]
For the Greenwood quantities themselves, positivity at finite sample size is
not automatic.  At fixed time points for which \(h(t)>0\), however, the
consistency results below imply that the corresponding scaled Greenwood
component is positive with probability tending to one; under the maintained
risk-set positivity assumptions the same applies to the scaled Greenwood
estimator.  Thus the log-scale calculations are asymptotically well defined in
the regime \(h(t)>0\), while deterministic scaling by \(n\) leaves their
log-scale variance unchanged.

The main technical point is that the first-order behavior of
\(\widetilde H_n(t)\) is not exhausted by the martingale fluctuation of the
event counting process.  The empirical at-risk process appearing in the
denominator is itself random, and its fluctuation contributes at the same
\(n^{-1/2}\) order.  More precisely, Theorem~\ref{thm:CLT_tildeH} yields the
asymptotically linear representation
\[
\sqrt n\{\widetilde H_n(t)-h(t)\}
=
\frac{1}{\sqrt n}\sum_{i=1}^n \phi_{H,t}(Z_i)+o_p(1),
\]
with the influence function \(\phi_{H,t}\) defined explicitly in that theorem.
It contains both the martingale event-process term and an additional contribution
generated by fluctuations of the empirical risk set.  The latter is part of the
full first-order asymptotic variance and cannot in general be omitted.
Accordingly, the focus here is not the consistency of \(\widetilde H_n(t)\)
itself, but its pointwise \(\sqrt n\)-first-order theory and its consequences for
Greenwood's variance estimator.

The influence-function representation also provides a direct route to the
asymptotic variance.  Let \(F_Z\) denote the distribution of
\(Z=(X,\Delta)\).  Since the influence function is centered, the variance
appearing in Theorem~\ref{thm:CLT_tildeH} is
\[
v_H(t)
=
\operatorname{Var}\{\phi_{H,t}(Z)\}
=
\mathbb E\{\phi_{H,t}(Z)^2\}
=
\int \phi_{H,t}(z)^2\,dF_Z(z),
\]
so that the first-order variance can be calculated directly from the squared
influence function.  The same representation is also the link from the
Greenwood component back to Greenwood's estimator.  A marginal central limit
theorem for \(\widetilde H_n(t)\) alone is not sufficient for a delta method
analysis of
\[
g(x,y)=x^2y,
\qquad
\widetilde G_n(t)=g\{\widehat S_n(t),\widetilde H_n(t)\},
\]
because the two arguments fluctuate jointly.  Combining the influence
function of \(\widetilde H_n(t)\) with the standard first-order representation
of \(\widehat S_n(t)\) \cite{reid1981,loSingh1986} yields the joint central limit theorem in
Theorem~\ref{thm:joint-AN-S-Htilde}, including the cross-covariance required
for the multivariate delta method.  This in turn gives the first-order
asymptotic distribution of the scaled Greenwood estimator; see
Corollary~\ref{cor:scaled-greenwood}.

Accordingly, the paper develops the pointwise first-order theory in the
following sequence.  We first establish the stochastic bounds and remainder
control needed for \(\widetilde H_n(t)\), record its consistency, and derive
its asymptotically linear representation, central limit theorem, and full
asymptotic variance.  We then combine this representation with the
Kaplan--Meier influence function to obtain joint asymptotic normality of
\((\widehat S_n(t),\widetilde H_n(t))\), and finally transfer the result to
\(\widetilde G_n(t)\) by the multivariate delta method.  The main theory is
presented under a continuous-time no-ties model.  The role of this assumption
and possible extensions to tied event times are discussed separately.
\section{Setup, Assumptions, and Greenwood Decomposition}
\subsection{Setup and Notation}

Fix a finite time horizon \(\tau>0\). Let \(T\ge0\) denote the event time with
survival function \(S(t)=\Pr(T>t)\), subject to right censoring by \(C\).
We assume \(C\ge0\) almost surely, allowing \(C=\infty\) for no censoring.
For a generic subject, set
\[
X=\min(T,C),
\qquad
\Delta=\mathbf{1}_{\{T\le C\}},
\qquad
Z=(X,\Delta),
\]
and let \(F_Z\) denote the distribution of \(Z\). For subjects
\(i=1,\dots,n\), define analogously
\[
X_i=\min(T_i,C_i),
\qquad
\Delta_i=\mathbf{1}_{\{T_i\le C_i\}},
\qquad
Z_i=(X_i,\Delta_i).
\]
The sampling and independent-censoring assumptions governing
\((T_i,C_i)\) are stated in Assumption~(A1) below.

For \(t\le \tau\), define the individual and aggregate at-risk processes by
\[
Y_i(t)=\mathbf{1}_{\{X_i\ge t\}},
\qquad
Y_n(t)=\sum_{i=1}^n Y_i(t).
\]
Thus \(Y_i(t)\) indicates whether subject \(i\) is still at risk at time \(t\),
while \(Y_n(t)\) is the corresponding number of subjects at risk.  Similarly,
define the individual and aggregate failure counting processes by
\[
N_i(t)=\mathbf{1}_{\{X_i\le t,\ \Delta_i=1\}},
\qquad
N_n(t)=\sum_{i=1}^n N_i(t).
\]
Here \(N_i(t)\) indicates whether an observed failure of subject \(i\) has
occurred by time \(t\), and \(N_n(t)\) is the total number of observed failures
by that time.

Let \(D_i^C(u):=\mathbf{1}_{\{X_i\le u,\ \Delta_i=0\}}\) denote the observed
censoring counting process. We take \(\mathcal F_t\) to be the usual completed,
right-continuous augmentation of the observed-history filtration
\[
\mathcal F_t^0
:=\sigma\{N_i(s),D_i^C(s):0\le s\le t,\ i=1,\dots,n\},
\qquad 0\le t\le\tau.
\]
Under the nonnegative-time convention, \(Y_i(0)=1\) and, for \(u>0\),
\[
Y_i(u)=1-N_i(u-)-D_i^C(u-).
\]
Thus \(Y_i\) and \(Y_n\) are left-continuous and adapted, hence predictable;
\(N_i\) and \(N_n\) are adapted, right-continuous counting processes. In
particular, the risk indicators are valid predictable integrands in the
counting-process martingale representations used below.

The corresponding population at-risk probability is denoted by
\[
y(t):=\Pr(X\ge t)=E\{Y(t)\},
\qquad
Y(t):=\mathbf{1}_{\{X\ge t\}}.
\]

Let \(0<t_1<t_2<\cdots\) denote the distinct observed failure times.  With
\[
\Delta N_n(t):=N_n(t)-N_n(t-),
\]
write
\[
d_j:=\Delta N_n(t_j),
\qquad
n_j:=Y_n(t_j-).
\]
Thus \(d_j\) is the number of observed failures at \(t_j\), whereas \(n_j\) is
the size of the risk set immediately before \(t_j\).  Since \(Y_n\) is
left-continuous, \(Y_n(t_j-)=Y_n(t_j)\); the notation \(t_j-\) is retained to
emphasize the usual pre-failure risk set.  Under Assumption~(A4) below,
\(d_j=1\) for every observed failure time.

The Kaplan--Meier product-limit estimator \cite{kaplan1958} is
\[
\widehat S_n(t)
=
\prod_{t_j\le t}
\left(1-\frac{d_j}{n_j}\right),
\qquad t\le \tau.
\]
Inference for \(\widehat S_n\) commonly employs Greenwood's variance estimator
\cite{greenwood1926}, which can be written as
\begin{equation}
\widehat G_n(t)
=
\widehat S_n(t)^2\,\widehat H_n(t),
\qquad t\le \tau,
\end{equation}
where the Greenwood component is defined pointwise in time as follows. For
\(t\le\tau\), let
\[
\mathcal D_n(t)
:=
\bigl\{n_j-d_j>0\ \text{for every observed failure time }t_j\le t\bigr\},
\]
and set
\[
\widehat H_n(t)
:=
\begin{cases}
\displaystyle\sum_{t_j\le t}\frac{d_j}{n_j(n_j-d_j)},
&\text{on }\mathcal D_n(t),\\[1.1em]
0,&\text{on }\mathcal D_n(t)^c.
\end{cases}
\]
The value zero on \(\mathcal D_n(t)^c\) is only a convenient arbitrary
extension; any fixed extension would be asymptotically equivalent.  Write
\(\mathcal D_n:=\mathcal D_n(\tau)\). Under Assumptions~(A1), (A2), and
(A4), the law of large numbers gives
\[
\frac{Y_n(\tau)}{n}\longrightarrow y(\tau)>0
\qquad\text{almost surely},
\]
so that \(\Pr\{Y_n(\tau)\ge2\}\to1\). Because \(Y_n\) is
nonincreasing and \(d_j=1\), on this event
\(n_j\ge Y_n(\tau)\ge2\) at every failure time \(t_j\le\tau\), and
therefore \(\{Y_n(\tau)\ge2\}\subseteq\mathcal D_n\). Hence
\(\Pr(\mathcal D_n)\to1\), and consequently
\(\Pr\{\mathcal D_n(t)\}\to1\) for every fixed \(t\le\tau\). All algebraic
identities involving Greenwood denominators up to time \(t\) below are
understood on \(\mathcal D_n(t)\); the arbitrary extension merely makes the
random quantities globally defined.

In the continuous-time no-ties model this becomes, on \(\mathcal D_n(t)\),
\[
\widehat H_n(t)
=
\sum_{t_j\le t}
\frac{1}{n_j(n_j-1)}.
\]

The first-order analysis below is therefore formulated in terms of the scaled
Greenwood component
\[
\widetilde H_n(t):=n\,\widehat H_n(t),
\]
and, correspondingly, the scaled Greenwood estimator
\[
\widetilde G_n(t)
:=
n\,\widehat G_n(t)
=
\widehat S_n(t)^2\,\widetilde H_n(t).
\]
The deterministic limit of \(\widetilde H_n(t)\) and its first-order
asymptotic properties are developed in Section~\ref{sec:first-order-H}.

\subsection{Regularity Conditions}
On the fixed finite horizon \([0,\tau]\) introduced above, the assumptions below
separate the sampling and censoring structure, positivity of the risk set,
regularity of the event-time hazard, and the continuous-time no-ties model.
\begin{enumerate}\setlength{\itemsep}{2pt}
\item[(A1)] \textbf{Sampling and independent censoring}
The pairs \((T_i,C_i)\), \(i=1,\dots,n\), are i.i.d., with \(T_i\perp C_i\).
Consequently, the observed data \(Z_i=(X_i,\Delta_i)\) are i.i.d.

\item[(A2)] \textbf{Positive at-risk probability:}
With \(y(t)=P(X\ge t)\) as defined above,
\[
\inf_{t\in[0,\tau]} y(t)>0.
\]
Thus the population at-risk probability stays bounded away from zero on
\([0,\tau]\). Since \(t\mapsto y(t)\) is nonincreasing, this condition is
equivalent to \(y(\tau)>0\). For a result at a single fixed time \(t\), the
corresponding local requirement \(y(t)>0\) would suffice; the horizon-wide
formulation is used to keep the assumptions uniform for all \(t\le\tau\).

\item[(A3)] \textbf{Event-hazard regularity:}
The event-time distribution is absolutely continuous on \([0,\tau]\), with
hazard function \(\lambda\). Together with this absolute continuity, (A1) and
(A2) imply \(S(\tau)>0\): by independent censoring, \(y(\tau)>0\) entails
\(P(T\ge\tau)>0\), and absolute continuity gives
\(P(T\ge\tau)=P(T>\tau)=S(\tau)\). Hence the cumulative hazard
\[
\Lambda(t):=\int_0^t \lambda(u)\,du
\]
is finite for every \(t\le\tau\), in particular
\(\Lambda(\tau)=-\log S(\tau)<\infty\). Boundedness of \(\lambda\) on
\([0,\tau]\) would of course be sufficient, since it implies
\(\Lambda(\tau)\le \tau\sup_{u\le\tau}\lambda(u)<\infty\), but it is a
strictly stronger condition than is needed here: an integrable hazard may be
unbounded, and the arguments below use only finiteness of the cumulative
hazard.

\item[(A4)] \textbf{Continuous-time no-ties model:}
Observed failure times are almost surely distinct, so that at each distinct
observed failure time \(t_j\le\tau\),
\[
d_j=\Delta N_n(t_j)=1.
\]
For exact continuous-time observations this property already follows from the
absolute continuity in (A3); it is retained as a separate condition to make the
no-ties structure used in the Greenwood decomposition explicit. Ties created
by rounding, grouping, or discrete recording of event times are not part of
the present asymptotic model. No continuity assumption on the censoring-time
distribution is required: (A1) and (A3) imply
\mbox{\(P(T=C\le\tau)=0\)}.
\end{enumerate}

The filtration introduced in Section~2.1 now gives the standard
multiplicative intensity representation explicitly. Under (A1)--(A3), the
observed failure process \(N_i\) has predictable intensity
\(Y_i(t)\lambda(t)\), and therefore
\[
M_i(t)
:=N_i(t)-\int_0^tY_i(u)\lambda(u)\,du
\]
is a square-integrable \((\mathcal F_t)\)-martingale with predictable
variation
\[
\langle M_i\rangle(t)
=\int_0^tY_i(u)\lambda(u)\,du.
\]
For a generic subject we use the analogous notation
\[
M(t):=N(t)-\int_0^tY(u)\lambda(u)\,du.
\]
Summing over subjects gives
\[
M_n(t):=\sum_{i=1}^n M_i(t),
\qquad
N_n(t)=M_n(t)+\int_0^tY_n(u)\lambda(u)\,du.
\]
The compensators are continuous, so that \(\Delta M_i=\Delta N_i\).
By (A1) and (A3), distinct subjects have no simultaneous observed
failure jumps on \([0,\tau]\). Thus, for \(i\ne k\),
\([M_i,M_k]\equiv0\) and hence \(\langle M_i,M_k\rangle\equiv0\).
Consequently,
\[
\langle M_n\rangle(t)=\int_0^tY_n(u)\lambda(u)\,du.
\]
Thus the risk process introduced in Section~2.1 is not only a bookkeeping
device: its predictability makes it the admissible pre-jump multiplier in the
counting-process martingales used below \cite{aalen1978,andersen1993,fleming1991}.
The assumptions needed for each subsequent result are stated explicitly.

\subsection{Greenwood's Variance Representation and Decomposition}
\label{subsec:notation-ties}
The Greenwood component admits an exact algebraic decomposition before the
no-ties assumption is used. On \(\mathcal D_n(t)\), for every observed failure
time \(t_j\le t\),
\[
\frac{d_j}{n_j(n_j-d_j)}
=
\frac{d_j}{n_j^2}
+
\frac{d_j^2}{n_j^2(n_j-d_j)}.
\]
Accordingly, define the scaled Nelson--Aalen plug-in variance term globally by
\[
U_n(t)
:=
n\sum_{t_j\le t}\frac{d_j}{n_j^2}.
\]
For the equivalent counting-process representation, since \(Y_n\) is
integer-valued, define the reciprocal risk-set weights on
\(\mathbb N_0=\{0,1,\ldots\}\) by
\[
w_1(x):=
\begin{cases}
x^{-1},&x\ge1,\\
0,&x=0,
\end{cases}
\qquad
w_2(x):=
\begin{cases}
x^{-2},&x\ge1,\\
0,&x=0,
\end{cases}
\]
and
\[
w_G(x):=
\begin{cases}
\{x(x-1)\}^{-1},&x\ge2,\\
0,&x\in\{0,1\}.
\end{cases}
\]
As Borel functions of the predictable process \(Y_n\), the processes
\(w_1\{Y_n(s)\}\), \(w_2\{Y_n(s)\}\), and \(w_G\{Y_n(s)\}\) are predictable. Hence
\begin{equation}\label{eq:Un-integral}
U_n(t)
=
n\int_0^t w_2\{Y_n(s)\}\,dN_n(s),
\end{equation}
which remains valid whether or not tied failures occur.

To make the Greenwood--Nelson--Aalen decomposition exact on the whole sample
space, define
\[
R_n(t)
:=
\widehat H_n(t)-\frac1n U_n(t).
\]
Then, identically,
\begin{equation}\label{eq:Htilde-decomposition}
\widetilde H_n(t)
=
U_n(t)+nR_n(t).
\end{equation}
On \(\mathcal D_n(t)\), the preceding algebra gives the explicit representation
\begin{equation}
R_n(t)
=
\sum_{t_j\le t}
\frac{d_j^2}{n_j^2(n_j-d_j)}.
\end{equation}
Thus the algebraic decomposition itself does not require the no-ties
assumption; no ties simplify the form and order of the remainder.

Under Assumption~\textup{(A4)}, \(d_j=1\) at every observed failure time. On
\(\mathcal D_n(t)\), this yields
\begin{equation}\label{eq:Hhat-discrete-ties}
\widehat H_n(t)
=
\sum_{t_j\le t}\frac{1}{n_j(n_j-1)},
\end{equation}
and
\[
R_n(t)
=
\sum_{t_j\le t}\frac{1}{n_j^2(n_j-1)}.
\]
Using the jump identity
\(\int_0^t f(s)\,dN_n(s)=\sum_{s\le t}f(s)\Delta N_n(s)\), the no-ties
Greenwood component also has the predictable integral representation
\begin{equation}
\widehat H_n(t)
=
\int_0^t w_G\{Y_n(s)\}\,dN_n(s),
\qquad\text{on }\mathcal D_n(t).
\end{equation}
The no-ties condition is technically useful precisely here. With tied failures,
the exact Greenwood denominator contains the contemporaneous jump size
\(d_j=\Delta N_n(t_j)\), which is not a pre-jump predictable quantity; under
(A4), that jump size is fixed at one, so the Greenwood weight depends only on
the predictable pre-jump risk set.

The first term in the decomposition,
\(\sum_{t_j\le t}d_j/n_j^2\), is the usual plug-in variance estimator for the
Nelson--Aalen estimator \cite{nelson1972,aalen1978,andersen1993}, and \(U_n\) is its scaled version. Thus \(R_n(t)\) is
the Greenwood--Nelson--Aalen remainder. It is not identically zero even under
the no-ties model. On the high-probability events developed in Section~3, the
risk sets satisfy \(n_j\asymp n\) uniformly. Each no-ties remainder summand is
then of order \(n^{-3}\), while at most \(n\) failures can occur, suggesting the
order \(R_n(t)=O_p(n^{-2})\); Lemma~\ref{lem:greenwood-remainder} makes this
heuristic precise.

This decomposition is the starting point for the first-order analysis in
Section~\ref{sec:first-order-H}. There we show that the remainder is negligible
at the first-order scale,
\[
\sqrt n\,\{nR_n(t)\}=o_p(1),
\]
so that, with \(h(t)\) denoting the deterministic limit introduced below,
\[
\sqrt n\{\widetilde H_n(t)-h(t)\}
=
\sqrt n\{U_n(t)-h(t)\}+o_p(1).
\]
Consequently, at each fixed \(t\), the first-order asymptotic behavior of
\(\widetilde H_n(t)\) is determined by that of \(U_n(t)\).

\section{First-Order Theory of the Greenwood Component}
\label{sec:first-order-H}

This section develops the first-order theory of the scaled Greenwood component
\(\widetilde H_n(t)\). The basic strategy is to use the exact decomposition
from Section~\ref{subsec:notation-ties} to replace \(\widetilde H_n(t)\), at
first order, by the scaled Nelson--Aalen variance term \(U_n(t)\), and then to
analyze \(U_n(t)\) by counting-process and empirical-process arguments. The
proof proceeds in four stages. We first establish, for each fixed \(t\),
uniform sample-level control on \([0,t]\) of the at-risk process and the
resulting Greenwood denominator bounds. We then show that the remainder in the
exact Greenwood--Nelson--Aalen decomposition is negligible not only for
consistency but also at the \(\sqrt n\)-scale. Next, for completeness and as
preparation for the first-order expansion, we recover the deterministic limit
\[
h(t)=\int_0^t \frac{\lambda(s)}{y(s)}\,ds
\]
through \(U_n(t)\). Finally, we refine this consistency result to the
\(\sqrt n\)-scale, derive an asymptotic linear representation, and obtain the
pointwise central limit theorem for \(\widetilde H_n(t)\).

\subsection{Preliminary Results}

We first establish the uniform sample-level control of the at-risk process
needed to bound the random denominators in the Greenwood terms and in the
subsequent stochastic expansions.

\begin{lemma}[Uniform convergence and positivity of the at-risk process]
\label{lem:uniform-positivity}
Suppose Assumption \textup{(A1)} holds. Then, for each fixed \(t\in[0,\tau]\),
\[
\sup_{0\le u\le t}\left|\frac{Y_n(u)}{n}-y(u)\right| \xrightarrow{a.s.} 0.
\]
If, in addition, Assumption \textup{(A2)} holds, then for each fixed
\(t \le \tau\) there exists a constant \(c_t \in (0,1)\) such that
\[
\inf_{0\le u \le t} \frac{Y_n(u)}{n} \ge c_t
\quad\text{with probability tending to one.}
\]
In particular, taking \(t=\tau\), there exists \(c\in(0,1)\) such that
\[
\inf_{0\le u \le \tau} \frac{Y_n(u)}{n} \ge c
\quad\text{with probability tending to one.}
\]
\end{lemma}

\begin{proof}
By Assumption~\textup{(A1)}, \(X_1,\dots,X_n\) are i.i.d.\ copies of \(X\).
Only this i.i.d.\ sampling part of (A1) is used here; independent censoring
enters the counting-process representation later. The empirical at-risk
function is
\[
\widehat y_n(u):=\frac{Y_n(u)}{n}
=\frac1n\sum_{i=1}^n\mathbf 1\{X_i\ge u\}.
\]
The class of half-line indicators
\(\{x\mapsto\mathbf 1\{x\ge u\}:0\le u\le t\}\) is a VC class and hence
Glivenko--Cantelli \cite[Ch.~19]{vdvaart1998}. Therefore
\(\widehat y_n\) converges uniformly almost surely to
\(y(u)=P(X\ge u)\) on \([0,t]\). No continuity of the distribution of \(X\)
is required; in particular, possible point masses induced by censoring do not
affect this conclusion. This proves the first assertion.

Assume also \textup{(A2)}. Since the indicators
\(\mathbf 1\{X_i\ge u\}\) are nonincreasing in \(u\), we have, for each
fixed \(t\in[0,\tau]\),
\[
\inf_{0\le u\le t}\frac{Y_n(u)}n
=\frac{Y_n(t)}n
=\frac1n\sum_{i=1}^n\mathbf 1\{X_i\ge t\}
\xrightarrow{a.s.}y(t),
\]
where the convergence follows from the ordinary strong law of large numbers.
By Assumption~\textup{(A2)}, \(y(t)>0\); taking
\(c_t:=y(t)/2\in(0,1)\) therefore yields
\[
P\!\left\{\inf_{0\le u\le t}\frac{Y_n(u)}n\ge c_t\right\}
\longrightarrow 1.
\]
Taking \(t=\tau\) gives the final assertion with \(c=y(\tau)/2\).
\end{proof}

Thus the uniform convergence and positivity bound provide the sample-level
denominator control used throughout the arguments below. The consequences
needed for the Greenwood terms are recorded explicitly next.

\begin{corollary}[Consequences of uniform positivity]
\label{cor:uniform-positivity-consequences}
Suppose Assumptions \textup{(A1)}, \textup{(A2)}, and \textup{(A4)} hold.
Then, for each fixed \(t\le\tau\), there exists \(c>0\) such that, with
probability tending to one,
\begin{enumerate}\setlength{\itemsep}{2pt}
\item[(U1)] \textbf{Risk-set lower bound:}
\[
Y_n(u)\ge cn\quad\text{for all }0\le u\le t,
\]
and hence \(n_j=Y_n(t_j-)\ge cn\) at every observed failure time \(t_j\le t\).

\item[(U2)] \textbf{Greenwood denominator control:}
under \textup{(A4)}, \(d_j=1\), and for every observed failure time
\(t_j\le t\),
\[
\frac{d_j}{n_j}=\frac1{n_j}\le\frac1{cn}=O(n^{-1}).
\]
For all sufficiently large \(n\), the corresponding Greenwood denominators
satisfy uniformly
\[
\frac12c^2n^2
\le n_j(n_j-d_j)\le n^2.
\]
\end{enumerate}
\end{corollary}

\begin{proof}
Item \textup{(U1)} is Lemma~\ref{lem:uniform-positivity}. On the event in
\textup{(U1)}, and under \textup{(A4)}, every observed failure time
\(t_j\le t\) satisfies \(d_j=1\) and \(n_j\ge cn\). This gives the displayed
bound for \(d_j/n_j\). For all sufficiently large \(n\),
\(1-d_j/n_j\ge1/2\) at every such failure time, and hence
\[
n_j(n_j-d_j)
=n_j^2\left(1-\frac{d_j}{n_j}\right)
\ge\frac12c^2n^2,
\]
while \(n_j(n_j-d_j)\le n_j^2\le n^2\). The statements are vacuous, and
therefore still valid, if no failure occurs by time \(t\).
\end{proof}

For all sufficiently large \(n\), the event in \textup{(U1)} is contained in
\(\mathcal D_n(t)\). On the same event, the globally defined reciprocal weights
from Section~\ref{subsec:notation-ties} reduce on \([0,t]\) to their ordinary
forms,
\[
\begin{aligned}
w_1\{Y_n(u)\}&=Y_n(u)^{-1},\\
w_2\{Y_n(u)\}&=Y_n(u)^{-2},\\
w_G\{Y_n(u)\}&=\{Y_n(u)(Y_n(u)-1)\}^{-1}.
\end{aligned}
\]
Thus the global definitions used to cover rare zero-denominator samples agree,
with probability tending to one, with the usual Greenwood and Nelson--Aalen
expressions on the interval relevant for a fixed-time result.

The same risk-set lower bound will reappear in the consistency and first-order
arguments below, where it controls inverse powers of \(Y_n\).

\begin{remark}[Pointwise order of the Greenwood component]
The preceding denominator bounds also recover directly the stochastic order of
the unscaled Greenwood component. For all sufficiently large \(n\), the
high-probability event from
Corollary~\ref{cor:uniform-positivity-consequences} is contained in
\(\mathcal D_n(t)\), and therefore
\[
\widehat H_n(t)
=\sum_{t_j\le t}\frac{1}{n_j(n_j-1)}
\le \frac{C}{n^2}N_n(t)
\le \frac{C}{n}
\]
for a finite constant \(C\). Hence
\[
\widehat H_n(t)=O_p(n^{-1}).
\]
This is an upper stochastic-order statement only; it does not imply that
\(n\widehat H_n(t)\) converges to a positive limit. In particular, the limit
\(h(t)\) identified below may equal zero.
\end{remark}

\subsection{Remainder Control}

The exact decomposition in Section~2.3 reduces the first-order analysis of
\(\widetilde H_n\) to that of \(U_n\), provided that the remainder \(nR_n\)
is negligible at the relevant scale. The following lemma gives a bound strong
enough for both consistency and the subsequent \(\sqrt n\)-analysis.

\begin{lemma}[Negligibility of the Greenwood remainder]
\label{lem:greenwood-remainder}
For each fixed \(t\le\tau\), under Assumptions \textup{(A1)}, \textup{(A2)},
and \textup{(A4)},
\[
R_n(t)=O_p(n^{-2}),
\qquad
nR_n(t)=O_p(n^{-1}),
\]
and consequently
\[
\sqrt n\,\{nR_n(t)\}=o_p(1).
\]
Equivalently,
\[
\sqrt n\{\widetilde H_n(t)-U_n(t)\}=o_p(1).
\]
\end{lemma}

\begin{proof}
The no-ties assumption enters this proof only through the representation of the
remainder with \(d_j=1\); the risk-set bound itself follows from
\textup{(A1)}--\textup{(A2)}. As in Lemma~\ref{lem:uniform-positivity}, only the
i.i.d.\ sampling implication of \textup{(A1)}, not its independent-censoring
part, is used here.

By Corollary~\ref{cor:uniform-positivity-consequences}, there exists \(c>0\)
such that
\[
E_n(t):=\left\{\inf_{0\le u\le t}\frac{Y_n(u)}n\ge c\right\}
\]
has probability tending to one. On \(E_n(t)\), \(n_j\ge cn\) for all failure
times \(t_j\le t\). For sufficiently large \(n\), this implies
\(E_n(t)\subseteq\mathcal D_n(t)\), so the explicit no-ties representation of
\(R_n(t)\) from Section~\ref{subsec:notation-ties} applies. Moreover,
\(n_j-1\ge n_j/2\), and therefore
\[
\frac1{n_j^2(n_j-1)}
\le\frac2{n_j^3}
\le\frac{2}{c^3n^3}.
\]
Since each subject can fail at most once, \(N_n(t)\le n\). Hence, on
\(E_n(t)\),
\[
R_n(t)
\le\frac{2N_n(t)}{c^3n^3}
\le\frac{2}{c^3n^2}.
\]
Because the explicit remainder representation is used only on \(E_n(t)\), the
arbitrary extension on its rare complement does not affect the stochastic
order: for all sufficiently large \(n\),
\[
P\left\{n^2|R_n(t)|>\frac{2}{c^3}\right\}
\le P\{E_n(t)^c\}\longrightarrow0.
\]
Thus \(R_n(t)=O_p(n^{-2})\). Consequently,
\[
nR_n(t)=O_p(n^{-1})=o_p(1),
\qquad
\sqrt n\,\{nR_n(t)\}=O_p(n^{-1/2})=o_p(1).
\]
Finally, \(\widetilde H_n(t)-U_n(t)=nR_n(t)\) from
\eqref{eq:Htilde-decomposition} gives the equivalent assertion in the lemma.
\end{proof}

Thus the same exact Greenwood--Nelson--Aalen decomposition serves both scales:
\(nR_n=o_p(1)\) is sufficient for consistency, while
\(\sqrt n\,\{nR_n\}=o_p(1)\) allows the first-order expansion of \(U_n\) to be
transferred unchanged to \(\widetilde H_n\).

\subsection{Consistency}
\label{subsec:consist}

Consistency of the scaled Greenwood component is already available from the
existing large-sample theory; in particular, it is contained in Gill's
product-limit theory \cite[Theorem~4.2.2]{gill1980}. For completeness, and to
prepare the first-order expansion below, we establish the pointwise result
directly from the Greenwood--Nelson--Aalen decomposition.

\begin{theorem}[Consistency of the scaled Greenwood component]
\label{thm:consistency-tildeH}
Let \(0\le t\le\tau\) be fixed. Under Assumptions \textup{(A1)}--\textup{(A4)},
\[
\widetilde H_n(t)\xrightarrow{P}h(t),
\qquad
h(t):=\int_0^t\frac{\lambda(s)}{y(s)}\,ds.
\]
Moreover, \(0\le h(t)<\infty\), and
\[
h(t)>0
\quad\Longleftrightarrow\quad
\Lambda(t)=\int_0^t\lambda(s)\,ds>0.
\]
\end{theorem}

\begin{proof}
First consider the scaled Nelson--Aalen variance term
\[
U_n(t)=n\int_0^t w_2\{Y_n(s)\}\,dN_n(s),
\]
where the globally defined predictable weight \(w_2\) is from
Section~\ref{subsec:notation-ties}. Using the multiplicative-intensity
representation established in Section~2.2 gives
\[
U_n(t)=A_n(t)+B_n(t),
\]
where, using \(Y_n w_2(Y_n)=w_1(Y_n)\),
\[
A_n(t):=n\int_0^t w_2\{Y_n(s)\}\,dM_n(s),
\qquad
B_n(t):=n\int_0^t w_1\{Y_n(s)\}\lambda(s)\,ds.
\]
These expressions are defined on the whole sample space, including rare samples
on which the risk set reaches zero.

Let
\[
m_t:=\inf_{0\le s\le t}y(s)>0,
\qquad
c:=m_t/2,
\qquad
E_n(t):=\left\{\inf_{0\le s\le t}\frac{Y_n(s)}n\ge c\right\}.
\]
By Lemma~\ref{lem:uniform-positivity}, \(P\{E_n(t)\}\to1\). To control the
martingale term, introduce the localized integral
\[
\begin{aligned}
A_n^{(c)}(t)
&:=n\int_0^t w_2\{Y_n(s)\}\mathbf 1\{Y_n(s)\ge cn\}\,dM_n(s)\\
&=\int_0^t\frac{n}{Y_n(s)^2}\mathbf 1\{Y_n(s)\ge cn\}\,dM_n(s).
\end{aligned}
\]
Its predictable variation satisfies
\[
\left\langle A_n^{(c)}\right\rangle(t)
=
\int_0^t
\frac{n^2}{Y_n(s)^3}\mathbf 1\{Y_n(s)\ge cn\}\lambda(s)\,ds
\le \frac{\Lambda(t)}{c^3n}.
\]
For later use, recall the martingale \(L^2\)-isometry: if \(M\) is a
square-integrable martingale and \(H\) is predictable with
\(E\int_0^t H(s)^2\,d\langle M\rangle(s)<\infty\), then
\[
E\left\{\left(\int_0^t H(s)\,dM(s)\right)^2\right\}
=
E\left\{\int_0^t H(s)^2\,d\langle M\rangle(s)\right\}.
\]
This is the square-integrable-martingale analogue of It\^o's isometry
\cite{andersen1993,fleming1991}. Hence the martingale isometry gives
\(E\{A_n^{(c)}(t)^2\}=O(n^{-1})\), and therefore
\(A_n^{(c)}(t)=o_p(1)\). Since \(A_n^{(c)}(t)=A_n(t)\) on \(E_n(t)\), for every
\(\varepsilon>0\),
\[
P\{|A_n(t)|>\varepsilon\}
\le
P\{E_n(t)^c\}+P\{|A_n^{(c)}(t)|>\varepsilon\}
\longrightarrow0.
\]
Thus \(A_n(t)=o_p(1)\).

For the drift term, write \(y_n(s):=Y_n(s)/n\). On \(E_n(t)\),
\(y_n(s)\ge m_t/2\) for all \(0\le s\le t\), so
\(n w_1\{Y_n(s)\}=1/y_n(s)\) and
\[
\begin{aligned}
\sup_{0\le s\le t}
\left|n\,w_1\{Y_n(s)\}-\frac1{y(s)}\right|
&=
\sup_{0\le s\le t}
\frac{|y_n(s)-y(s)|}{y_n(s)y(s)}\\
&\le
\frac{2}{m_t^2}
\sup_{0\le s\le t}|y_n(s)-y(s)|.
\end{aligned}
\]
Lemma~\ref{lem:uniform-positivity} and \(P\{E_n(t)^c\}\to0\) therefore imply
\[
\sup_{0\le s\le t}
\left|n\,w_1\{Y_n(s)\}-\frac1{y(s)}\right|
\xrightarrow{P}0.
\]
Consequently,
\[
\begin{aligned}
|B_n(t)-h(t)|
&\le
\Lambda(t)
\sup_{0\le s\le t}
\left|n\,w_1\{Y_n(s)\}-\frac1{y(s)}\right|
\xrightarrow{P}0.
\end{aligned}
\]
Thus \(U_n(t)\xrightarrow{P}h(t)\). This part of the argument uses only
Assumptions \textup{(A1)}--\textup{(A3)}. Assumption \textup{(A4)} enters through
Lemma~\ref{lem:greenwood-remainder}, which gives
\[
\widetilde H_n(t)=U_n(t)+o_p(1),
\]
and hence \(\widetilde H_n(t)\xrightarrow{P}h(t)\).

Finally, since \(0<y(s)\le1\) on \([0,t]\),
\[
0\le \Lambda(t)
\le h(t)
\le \frac{\Lambda(t)}{m_t}
<\infty.
\]
Therefore \(h(t)>0\) if and only if \(\Lambda(t)>0\), equivalently
\(h(t)=0\) if and only if \(\Lambda(t)=0\).
\end{proof}

The case \(h(t)=0\) therefore causes no difficulty for the additive first-order
theory. In that case \(\Lambda(t)=0\), and under (A3),
\(S(t)=\exp\{-\Lambda(t)\}=1\), so \(P(T\le t)=0\). Hence no failure occurs by
time \(t\) with probability one and, for every \(n\),
\[
N_n(t)=U_n(t)=\widehat H_n(t)=\widetilde H_n(t)=0
\qquad\text{almost surely}.
\]
The influence function in Theorem~\ref{thm:CLT_tildeH} is then zero almost
surely, so the pointwise central limit theorem is valid in its degenerate form
with \(v_H(t)=0\). Positivity is required only for transformations such as the
logarithm; the corresponding statement in Section~4.2 is therefore explicitly
restricted to \(h(t)>0\).

\subsection{Asymptotic Linearity and Central Limit Theorem}

Theorem~\ref{thm:consistency-tildeH} established the deterministic first-order
limit through \(U_n(t)\xrightarrow{P}h(t)\). We now refine this result to the
\(\sqrt n\)-scale. The proof is organized as follows. First, a Taylor expansion
of the functional representation of \(U_n-h\) separates the two linear
empirical fluctuations that contribute to first order from two genuinely
second-order terms. We then control the Taylor remainder directly and the
empirical cross term through a technical lemma. This yields the influence
function of \(U_n\). Finally, Lemma~\ref{lem:greenwood-remainder} transfers the
same expansion to \(\widetilde H_n\), after which the ordinary i.i.d.\ central
limit theorem applies.

\begin{lemma}[Empirical cross term]
\label{lem:empirical-cross-term}
Under Assumptions \textup{(A1)}--\textup{(A3)}, for each fixed \(t\le\tau\),
with
\[
y_n(s):=\frac{Y_n(s)}n,
\qquad
\nu_n(ds):=\frac1n\,dN_n(s),
\qquad
\nu(ds):=y(s)\lambda(s)\,ds,
\]
we have
\[
\int_0^t
\frac{y_n(s)-y(s)}{y(s)^3}\{\nu_n-\nu\}(ds)
=O_p(n^{-1}).
\]
\end{lemma}

The lemma shows that the simultaneous empirical fluctuations of the at-risk
process and the event measure interact only at second order. Consequently,
this product term does not contribute to the influence function. A detailed
second-moment proof is given in Appendix~\ref{app:cross-term-proof}.

\begin{theorem}[Asymptotic linearity and pointwise CLT for the scaled Greenwood component]
\label{thm:CLT_tildeH}
Suppose Assumptions \textup{(A1)}--\textup{(A4)} hold, and let
\(t\in[0,\tau]\) be fixed. For the generic observation \(Z=(X,\Delta)\), let
\(Y\), \(N\), and \(M\) be the processes defined in Sections~2.1--2.2. Set
\begin{equation}\label{eq:phi-H-theorem}
\phi_{H,t}(Z)
:=
\int_0^t\frac1{y(u)^2}\,dM(u)
-
\int_0^t\frac{\lambda(u)}{y(u)^2}\{Y(u)-y(u)\}\,du.
\end{equation}
Then \(E\{\phi_{H,t}(Z)\}=0\), \(E\{\phi_{H,t}(Z)^2\}<\infty\), and
\begin{equation}\label{eq:H-asymptotic-linearity}
\sqrt n\{\widetilde H_n(t)-h(t)\}
=
\frac1{\sqrt n}\sum_{i=1}^n\phi_{H,t}(Z_i)+o_p(1).
\end{equation}
Consequently,
\[
\sqrt n\{\widetilde H_n(t)-h(t)\}
\xRightarrow{d}\mathcal N\{0,v_H(t)\},
\]
where
\begin{equation}
v_H(t)
:=\operatorname{Var}\{\phi_{H,t}(Z)\}
=E\{\phi_{H,t}(Z)^2\}.
\end{equation}
\end{theorem}

\begin{proof}
Write
\[
y_n(s):=\frac{Y_n(s)}n,
\qquad
\nu_n(ds):=\frac1n\,dN_n(s),
\qquad
\nu(ds):=y(s)\lambda(s)\,ds.
\]
Let
\[
m_t:=\inf_{0\le s\le t}y(s)>0.
\]
The same VC class of half-line indicators used in
Lemma~\ref{lem:uniform-positivity} is Donsker
\cite{vdvaartWellner1996}; see also \cite[Ch.~19]{vdvaart1998}. Thus,
\[
\sup_{0\le s\le t}|y_n(s)-y(s)|=O_p(n^{-1/2}).
\]
Define the high-probability event
\[
\mathcal E_n(t):=
\left\{\sup_{0\le s\le t}|y_n(s)-y(s)|<m_t/2\right\}.
\]
Then \(P\{\mathcal E_n(t)\}\to1\), and on \(\mathcal E_n(t)\),
\(y_n(s)\ge m_t/2>0\) uniformly for \(0\le s\le t\). Hence the globally
defined representation \eqref{eq:Un-integral} reduces on this event to
\[
U_n(t)=\int_0^t\frac1{y_n(s)^2}\,\nu_n(ds),
\qquad
h(t)=\int_0^t\frac1{y(s)^2}\,\nu(ds).
\]
Thus the ordinary reciprocal representation is used only on
\(\mathcal E_n(t)\); outside this event, \(U_n\) remains defined through the
global weight \(w_2\).

On \(\mathcal E_n(t)\), Taylor expansion of \(f(x)=x^{-2}\) gives, uniformly
for \(0\le s\le t\),
\[
\frac1{y_n(s)^2}
=
\frac1{y(s)^2}
-2\frac{y_n(s)-y(s)}{y(s)^3}
+\rho_n(s).
\]
Since \(f''(x)=6x^{-4}\) and every intermediate value in Taylor's formula is
at least \(m_t/2\) on \(\mathcal E_n(t)\),
\[
|\rho_n(s)|
\le
48m_t^{-4}|y_n(s)-y(s)|^2.
\]
For definiteness, set \(\rho_n(s)=0\) on \(\mathcal E_n(t)^c\). Consequently,
\[
\sup_{0\le s\le t}|\rho_n(s)|
\le
48m_t^{-4}
\sup_{0\le s\le t}|y_n(s)-y(s)|^2
=
O_p(n^{-1}).
\]

Substitution into the representation of \(U_n\) on \(\mathcal E_n(t)\),
followed by adding and subtracting the corresponding \(\nu\)-term, gives
\begin{equation}\label{eq:key-Un-decomposition}
\begin{aligned}
U_n(t)-h(t)
={}&
\int_0^t\frac1{y(s)^2}\{\nu_n-\nu\}(ds)
-2\int_0^t\frac{y_n(s)-y(s)}{y(s)^3}\,\nu(ds)
\\
&\quad
-2\int_0^t\frac{y_n(s)-y(s)}{y(s)^3}\{\nu_n-\nu\}(ds)
+\int_0^t\rho_n(s)\,\nu_n(ds).
\end{aligned}
\end{equation}
This decomposition is the key step in the first-order analysis. The first term
is the linear fluctuation of the empirical event measure, and the second is the
linear fluctuation induced by the random at-risk process. These two terms
produce the two components of the influence function. The third term is their
empirical product and is therefore second order; the fourth is the Taylor
remainder.

For the Taylor remainder, \(\nu_n([0,t])=N_n(t)/n\le1\), and hence
\[
\left|\int_0^t\rho_n(s)\,\nu_n(ds)\right|
\le\sup_{0\le s\le t}|\rho_n(s)|
=O_p(n^{-1})=o_p(n^{-1/2}).
\]
By Lemma~\ref{lem:empirical-cross-term}, the third term in
\eqref{eq:key-Un-decomposition} is also \(O_p(n^{-1})\). Define the globally
well-defined linearization error
\[
\begin{aligned}
r_n(t):={}&U_n(t)-h(t)
-\int_0^t\frac1{y(s)^2}\{\nu_n-\nu\}(ds)\\
&\quad
+2\int_0^t\frac{\lambda(s)}{y(s)^2}\{y_n(s)-y(s)\}\,ds.
\end{aligned}
\]
On \(\mathcal E_n(t)\), the preceding decomposition gives
\[
r_n(t)
=
-2\int_0^t\frac{y_n(s)-y(s)}{y(s)^3}\{\nu_n-\nu\}(ds)
+\int_0^t\rho_n(s)\,\nu_n(ds).
\]
Thus, by Lemma~\ref{lem:empirical-cross-term} and the Taylor-remainder bound,
\(r_n(t)=O_p(n^{-1})\) on \(\mathcal E_n(t)\). Since
\(P\{\mathcal E_n(t)^c\}\to0\), for every \(\varepsilon>0\),
\[
\begin{aligned}
P\{\sqrt n\,|r_n(t)|>\varepsilon\}
&\le P\{\mathcal E_n(t)^c\}\\
&\quad+P\{2\sqrt n\,|C_n(t)|
+\sqrt n\sup_{0\le s\le t}|\rho_n(s)|>\varepsilon\}
\longrightarrow0,
\end{aligned}
\]
where \(C_n(t)\) denotes the empirical cross term in
Lemma~\ref{lem:empirical-cross-term}.
Therefore, on the whole sample space,
\[
U_n(t)-h(t)
=
\int_0^t\frac1{y(s)^2}\{\nu_n-\nu\}(ds)
-2\int_0^t\frac{\lambda(s)}{y(s)^2}\{y_n(s)-y(s)\}\,ds
+o_p(n^{-1/2}).
\]
Using
\[
\nu_n(ds)=\frac1n\sum_{i=1}^n dN_i(s),
\qquad
y_n(s)=\frac1n\sum_{i=1}^nY_i(s),
\]
we obtain
\begin{equation}\label{eq:Un-IF-expansion}
\sqrt n\{U_n(t)-h(t)\}
=
\frac1{\sqrt n}\sum_{i=1}^n\phi_{H,t}(Z_i)+o_p(1),
\end{equation}
where, equivalently,
\begin{align}
\phi_{H,t}(Z_i)
&=
\int_0^t\frac1{y(s)^2}\,dN_i(s)
-2\int_0^t\frac{\lambda(s)}{y(s)^2}Y_i(s)\,ds
+\int_0^t\frac{\lambda(s)}{y(s)}\,ds
\\
&=
\int_0^t\frac1{y(s)^2}\,dM_i(s)
-
\int_0^t\frac{\lambda(s)}{y(s)^2}\{Y_i(s)-y(s)\}\,ds.
\end{align}
Thus the complete influence function contains both the event-process
contribution and the random-risk-set contribution; the latter can
contribute at the \(\sqrt n\)-scale and cannot in general be omitted.

The martingale term has mean zero, and
\(E\{Y(s)-y(s)\}=0\), so
\[
E\{\phi_{H,t}(Z)\}=0.
\]
Moreover, let \(m_t:=\inf_{0\le u\le t}y(u)>0\). By the martingale isometry,
\[
E\left\{\left(\int_0^t\frac{1}{y(u)^2}\,dM(u)\right)^2\right\}
=\int_0^t\frac{\lambda(u)}{y(u)^3}\,du
\le m_t^{-3}\Lambda(t)<\infty.
\]
For the risk-set term, \(|Y(u)-y(u)|\le1\) gives
\[
\left|\int_0^t
\frac{\lambda(u)}{y(u)^2}\{Y(u)-y(u)\}\,du\right|
\le m_t^{-2}\Lambda(t)<\infty.
\]
Hence, by \((a-b)^2\le2a^2+2b^2\),
\(E\{\phi_{H,t}(Z)^2\}<\infty\).

The first-order expansion of \(U_n(t)\) itself uses only
Assumptions~\textup{(A1)}--\textup{(A3)}; Assumption~\textup{(A4)} enters only
through the Greenwood remainder bound used to transfer that expansion to
\(\widetilde H_n(t)\). Finally, Lemma~\ref{lem:greenwood-remainder} gives
\[
\sqrt n\{\widetilde H_n(t)-U_n(t)\}=o_p(1).
\]
Combining this with \eqref{eq:Un-IF-expansion} proves
\eqref{eq:H-asymptotic-linearity}. The ordinary i.i.d.\ central limit theorem
then yields
\[
\sqrt n\{\widetilde H_n(t)-h(t)\}
\xRightarrow{d}
\mathcal N\left(0,E\{\phi_{H,t}(Z)^2\}\right),
\]
which proves the theorem.
\end{proof}

\begin{remark}[Equivalent expanded-variance formulation]
The influence-function form in Theorem~\ref{thm:CLT_tildeH} makes explicit that
both the event-counting process and the empirical at-risk process enter the
influence-function representation and may contribute to the first-order
fluctuation of \(\widetilde H_n(t)\). An equivalent formulation
with the variance written out in integral form is given in
Appendix~\ref{app:expanded-H-clt}.
\end{remark}

\section{From the Greenwood Component to Greenwood's Estimator}
\subsection{Pointwise Joint Asymptotic Normality of
\texorpdfstring{$(\widehat S_n(t),\widetilde H_n(t))$}{(S-hat, H-tilde)}}
\label{subsec:joint-AN}

We now show that the Kaplan--Meier estimator and the scaled Greenwood
component are jointly asymptotically normal at each fixed time point.

\begin{theorem}[Joint CLT for
\texorpdfstring{$(\widehat S_n(t),\widetilde H_n(t))$}{(S-hat, H-tilde)}]
\label{thm:joint-AN-S-Htilde}
Suppose Assumptions \textup{(A1)}, \textup{(A2)}, \textup{(A3)}, and
\textup{(A4)} hold. Let \(t\in[0,\tau]\) be fixed. Let
\(\widehat S_n(t)\) denote the Kaplan--Meier estimator of \(S(t)\), and let
\[
\widetilde H_n(t):=n\,\widehat H_n(t),
\]
where \(\widehat H_n(t)\) is the globally defined Greenwood component from
Section~2.1, with the no-ties representation
\eqref{eq:Hhat-discrete-ties} on \(\mathcal D_n(t)\). Let \(Z=(X,\Delta)\),
and let \(Y\), \(N\), and \(M\) be the generic processes defined in
Sections~2.1--2.2. Define
\[
\psi_{S,t}(Z)
:=
-\,S(t)\int_0^t \frac{1}{y(u)}\,dM(u),
\]
and let \(\phi_{H,t}(Z)\) be the influence function in
\eqref{eq:phi-H-theorem}. Put
\[
v_S(t):=\operatorname{Var}\{\psi_{S,t}(Z)\}=S(t)^2h(t),
\qquad
c_{SH}(t):=\operatorname{Cov}\{\psi_{S,t}(Z),\phi_{H,t}(Z)\},
\]
and set
\[
\Sigma(t)
:=
\begin{pmatrix}
v_S(t) & c_{SH}(t)\\[2pt]
c_{SH}(t) & v_H(t)
\end{pmatrix}
=
\operatorname{Var}
\begin{pmatrix}
\psi_{S,t}(Z)\\[2pt]
\phi_{H,t}(Z)
\end{pmatrix}.
\]
Then
\[
\sqrt{n}
\begin{pmatrix}
\widehat S_n(t)-S(t)\\[2pt]
\widetilde H_n(t)-h(t)
\end{pmatrix}
\xRightarrow{d}
\mathcal N_2\bigl(0,\Sigma(t)\bigr),
\]
where \(h(t)\) is the deterministic limit of \(\widetilde H_n(t)\) from
Section~\ref{subsec:consist}. In particular, the lower-right entry of
\(\Sigma(t)\) is
\[
\operatorname{Var}\{\phi_{H,t}(Z)\}=v_H(t).
\]
No nonsingularity of \(\Sigma(t)\) is required. In particular, when
\(h(t)=0\), Section~\ref{subsec:consist} implies \(\Sigma(t)=0\), so the
joint limit is degenerate at the origin.
\end{theorem}

\begin{proof}
Under Assumptions~\textup{(A1)}--\textup{(A3)}, standard counting-process
theory for the Kaplan--Meier estimator gives the asymptotic linear
representation \cite{Wellner1982,gill1983,andersen1993,fleming1991}
\[
\sqrt{n}\bigl\{\widehat S_n(t)-S(t)\bigr\}
=
\frac{1}{\sqrt{n}}\sum_{i=1}^n \psi_{S,t}(Z_i)
+ o_p(1).
\]
Moreover, the martingale integral has mean zero, so
\[
E\{\psi_{S,t}(Z)\}=0.
\]
The martingale $L^2$-isometry introduced in
Section~\ref{subsec:consist} gives
\[
E\{\psi_{S,t}(Z)^2\}
=
S(t)^2\int_0^t\frac{\lambda(u)}{y(u)}\,du
=
S(t)^2h(t)<\infty.
\]
Thus the Kaplan--Meier expansion itself does not require the no-ties
assumption. Assumption~\textup{(A4)} enters the joint result through
Theorem~\ref{thm:CLT_tildeH}, which gives the corresponding asymptotic linear
representation for the scaled Greenwood component,
\[
\sqrt{n}\bigl\{\widetilde H_n(t)-h(t)\bigr\}
=
\frac{1}{\sqrt{n}}\sum_{i=1}^n \phi_{H,t}(Z_i)
+ o_p(1).
\]
Therefore,
\[
\sqrt{n}
\begin{pmatrix}
\widehat S_n(t)-S(t)\\[2pt]
\widetilde H_n(t)-h(t)
\end{pmatrix}
=
\frac{1}{\sqrt{n}}\sum_{i=1}^n
\begin{pmatrix}
\psi_{S,t}(Z_i)\\[2pt]
\phi_{H,t}(Z_i)
\end{pmatrix}
+o_p(1).
\]
By Theorem~\ref{thm:CLT_tildeH}, \(\phi_{H,t}(Z)\) is centered and square
integrable. Hence the two-dimensional summands above are i.i.d., centered, and
square integrable. The multivariate central limit theorem therefore gives
\[
\frac{1}{\sqrt{n}}\sum_{i=1}^n
\begin{pmatrix}
\psi_{S,t}(Z_i)\\[2pt]
\phi_{H,t}(Z_i)
\end{pmatrix}
\xRightarrow{d}
\mathcal N_2\bigl(0,\Sigma(t)\bigr),
\]
with
\[
\Sigma(t)
=
\operatorname{Var}
\begin{pmatrix}
\psi_{S,t}(Z)\\[2pt]
\phi_{H,t}(Z)
\end{pmatrix}.
\]
This proves the claim.
\end{proof}

\begin{remark}[Explicit cross-covariance]
The off-diagonal entry of \(\Sigma(t)\) can be written explicitly. Let
\[
J(t):=\int_0^t\frac{\lambda(u)}{y(u)^2}\,du.
\]
Since $m_t=\inf_{0\le u\le t}y(u)>0$,
\[
0\le J(t)\le m_t^{-2}\Lambda(t)<\infty.
\]
The martingale \(L^2\)-isometry introduced in Section~\ref{subsec:consist}
was stated there for a single integrand \(H\), and hence for the square of one
martingale integral. For two integrands \(H\) and \(K\), apply that quadratic
identity to \(H+K\) and \(H-K\) and subtract. Using the linearity of the
stochastic integral together with
\[
(H+K)^2-(H-K)^2=4HK,
\]
the cross-product terms are isolated, which gives
\[
E\left\{
\left(\int_0^t H(u)\,dM(u)\right)
\left(\int_0^t K(u)\,dM(u)\right)
\right\}
=
E\left\{\int_0^t H(u)K(u)\,d\langle M\rangle(u)\right\}.
\]
This passage from the quadratic identity to the corresponding cross-product,
or covariance, identity is called \emph{polarization}; the resulting formula
is the bilinear form of the martingale isometry. With \(H(u)=1/y(u)\) and
\(K(u)=1/y(u)^2\), it gives
\[
E\left\{
\left(\int_0^t\frac1{y(u)}\,dM(u)\right)
\left(\int_0^t\frac1{y(u)^2}\,dM(u)\right)
\right\}
=J(t).
\]
For the mixed moment
\(E\{Y(r)\int_0^t y(u)^{-1}\,dM(u)\}\), the event and drift parts are
\(\int_r^t\lambda(u)\,du\) and
\(y(r)h(r)+\int_r^t\lambda(u)\,du\), respectively, by the
atom-free measure \(\nu(du)=y(u)\lambda(u)\,du\) and
\(Y(r)Y(u)=Y(r\vee u)\). Thus
\[
E\left\{
Y(r)\int_0^t\frac1{y(u)}\,dM(u)
\right\}
=-y(r)h(r).
\]
The bounds \(\left|\int_0^t y^{-1}\,dM\right|
\le m_t^{-1}(1+\Lambda(t))\) and
\(\int_0^t\lambda/y^2\le m_t^{-2}\Lambda(t)\) justify
interchanging the expectation and the \(r\)-integral by Fubini;
the absolute continuity of \(h\) then gives
\[
E\left[
\left\{\int_0^t\frac1{y(u)}\,dM(u)\right\}
\left\{\int_0^t\frac{\lambda(r)}{y(r)^2}
      (Y(r)-y(r))\,dr\right\}
\right]
=-\frac12 h(t)^2,
\]
and hence
\[
\boxed{
 c_{SH}(t)
 =
 -S(t)\left\{
 J(t)+\frac12 h(t)^2
 \right\}.
}
\]
Thus the cross-covariance retained in the joint CLT is generally nonzero and
contributes to the delta-method variance of the scaled Greenwood estimator.
More precisely, \(c_{SH}(t)=0\) when \(h(t)=0\), whereas
\[
c_{SH}(t)<0
\qquad\text{whenever } h(t)>0.
\]
\end{remark}

\subsection{Delta-Method Consequences and the Scaled Greenwood Estimator}
\label{subsec:delta-method-consequences}

We briefly collect the conclusions of the preceding subsections. The uniform
positivity result ensures that the empirical risk set \(Y_n(u)\) remains of
order \(n\), uniformly on compact time intervals on which \(y(u)>0\). This
makes the denominators in Greenwood's formula asymptotically stable. The
remainder analysis then shows that the Greenwood component is asymptotically
equivalent, at the relevant scale, to the corresponding Nelson--Aalen plug-in
variance term. Consequently,
\[
\widetilde H_n(t)
=
n\,\widehat H_n(t)
\]
has the deterministic limit
\[
h(t)
=
\int_0^t \frac{\lambda(s)}{y(s)}\,ds,
\]
and Theorem~\ref{thm:CLT_tildeH} gives the pointwise central limit theorem
\[
\sqrt n\{\widetilde H_n(t)-h(t)\}
\xRightarrow{d}
\mathcal N\{0,v_H(t)\},
\qquad
v_H(t)
=
\operatorname{Var}\{\phi_{H,t}(Z)\}.
\]

If \(h(t)>0\), consistency implies
\[
P\{\widetilde H_n(t)\le0\}
\le
P\{|\widetilde H_n(t)-h(t)|\ge h(t)\}
\longrightarrow0,
\]
and hence
\[
P\{\widetilde H_n(t)>0\}\longrightarrow1.
\]
Thus the logarithm is well defined with probability tending to one; for
formal completeness it may be defined arbitrarily on the asymptotically
negligible complementary event. The ordinary one-dimensional delta method
\cite[Ch.~3]{vdvaart1998} then gives
\[
\sqrt n
\{\log \widetilde H_n(t)-\log h(t)\}
\xRightarrow{d}
\mathcal N\!\left(0,\frac{v_H(t)}{h(t)^2}\right).
\]
On the same positivity event,
\[
\log \widehat H_n(t)
=
\log \widetilde H_n(t)-\log n,
\]
so the deterministic shift \(-\log n\) affects the centering but not the
log-scale asymptotic variance. Equivalently,
\[
\sqrt n\left[
\log \widehat H_n(t)
-
\log\left\{\frac{h(t)}{n}\right\}
\right]
\xRightarrow{d}
\mathcal N\!\left(0,\frac{v_H(t)}{h(t)^2}\right).
\]

For Greenwood's full variance estimator, the relevant input is the joint CLT
\[
\sqrt n
\begin{pmatrix}
\widehat S_n(t)-S(t)\\[2pt]
\widetilde H_n(t)-h(t)
\end{pmatrix}
\xRightarrow{d}
\mathcal N_2\{0,\Sigma(t)\},
\]
from Theorem~\ref{thm:joint-AN-S-Htilde}. This is the precise input needed for
the multivariate delta method \cite[Ch.~3]{vdvaart1998} applied to the product map
\[
g(x,y)=x^2y.
\]
Thus the delta method is not applied to the unscaled pair
\((\widehat S_n(t),\widehat H_n(t))\), but to the correctly scaled pair
\((\widehat S_n(t),\widetilde H_n(t))\). This yields the asymptotic properties
of the scaled Greenwood estimator
\[
\widetilde G_n(t)
=
n\,\widehat G_n(t)
=
\widehat S_n(t)^2\,\widetilde H_n(t),
\]
as stated next.

\begin{corollary}[Asymptotic properties of the scaled Greenwood estimator]
\label{cor:scaled-greenwood}
Suppose Assumptions \textup{(A1)}, \textup{(A2)}, \textup{(A3)}, and
\textup{(A4)} hold. Let \(t\in[0,\tau]\) be fixed, and define
\[
\widetilde G_n(t)
:=
n\,\widehat G_n(t)
=
\widehat S_n(t)^2\,\widetilde H_n(t),
\qquad
\widetilde H_n(t):=n\,\widehat H_n(t).
\]
Let
\[
h(t)
=
\int_0^t \frac{\lambda(s)}{y(s)}\,ds
\]
denote the deterministic limit of \(\widetilde H_n(t)\), and set
\[
G_\ast(t):=S(t)^2h(t).
\]
Then
\[
\widetilde G_n(t)\xrightarrow{p}G_\ast(t).
\]
Moreover, define
\[
\phi_{G,t}(Z)
:=
2S(t)h(t)\psi_{S,t}(Z)
+
S(t)^2\phi_{H,t}(Z).
\]
Then \(E\{\phi_{G,t}(Z)\}=0\), \(E\{\phi_{G,t}(Z)^2\}<\infty\), and
\[
\sqrt n\{\widetilde G_n(t)-G_\ast(t)\}
=
\frac1{\sqrt n}\sum_{i=1}^n\phi_{G,t}(Z_i)+o_p(1)
\xRightarrow{d}
\mathcal N\bigl(0,v_G(t)\bigr),
\]
where
\[
v_G(t)
:=
\operatorname{Var}\{\phi_{G,t}(Z)\}
=
\dot g(t)^\top\Sigma(t)\dot g(t),
\qquad
\dot g(t)
=
\begin{pmatrix}
2S(t)h(t)\\[2pt]
S(t)^2
\end{pmatrix},
\]
and \(\Sigma(t)\) is the covariance matrix in
Theorem~\ref{thm:joint-AN-S-Htilde}. Consequently,
\[
\widehat G_n(t)
=
\frac{1}{n}\widetilde G_n(t)
=
O_p(n^{-1}).
\]
Thus the natural \(O_p(1)\)-scaled large-sample object is
\(\widetilde G_n(t)=n\widehat G_n(t)\), rather than \(\widehat G_n(t)\) itself.
\end{corollary}

\begin{proof}
By Theorem~\ref{thm:joint-AN-S-Htilde},
\[
\widehat S_n(t)\xrightarrow{p}S(t),
\qquad
\widetilde H_n(t)\xrightarrow{p}h(t).
\]
Since the map
\[
g(x,y)=x^2y
\]
is continuous, the continuous mapping theorem gives
\[
\widetilde G_n(t)
=
g\{\widehat S_n(t),\widetilde H_n(t)\}
\xrightarrow{p}
g\{S(t),h(t)\}
=
S(t)^2h(t)
=
G_\ast(t).
\]

The map \(g\) is continuously differentiable, with gradient
\[
\nabla g(x,y)
=
\begin{pmatrix}
2xy\\[2pt]
x^2
\end{pmatrix}.
\]
Evaluating this gradient at \(\{S(t),h(t)\}\) gives
\[
\nabla g\{S(t),h(t)\}
=
\begin{pmatrix}
2S(t)h(t)\\[2pt]
S(t)^2
\end{pmatrix}
=
\dot g(t).
\]
Applying the multivariate delta method to the asymptotically linear
representation in Theorem~\ref{thm:joint-AN-S-Htilde} yields the expansion
below. More precisely, the joint central limit theorem is stated in that
theorem, while the stacked asymptotically linear representation used here is
established in its proof; applying the first-order Taylor expansion of \(g\)
to that representation gives
\[
\begin{aligned}
\sqrt n
\left[
g\{\widehat S_n(t),\widetilde H_n(t)\}
-
g\{S(t),h(t)\}
\right]
&=
\dot g(t)^\top
\frac1{\sqrt n}\sum_{i=1}^n
\begin{pmatrix}
\psi_{S,t}(Z_i)\\[2pt]
\phi_{H,t}(Z_i)
\end{pmatrix}
+o_p(1)\\
&=
\frac1{\sqrt n}\sum_{i=1}^n\phi_{G,t}(Z_i)+o_p(1).
\end{aligned}
\]
Since
\[
g\{\widehat S_n(t),\widetilde H_n(t)\}
=
\widetilde G_n(t),
\qquad
g\{S(t),h(t)\}
=
G_\ast(t),
\]
the asserted asymptotic linearity follows. Since \(\psi_{S,t}(Z)\) and
\(\phi_{H,t}(Z)\) are centered and square integrable, the same holds for their
fixed linear combination \(\phi_{G,t}(Z)\). The ordinary i.i.d.\ central limit
theorem then gives the stated normal limit, with
\[
v_G(t)
=
\operatorname{Var}\{\phi_{G,t}(Z)\}
=
\dot g(t)^\top\Sigma(t)\dot g(t).
\]

Finally,
\[
\widehat G_n(t)
=
n^{-1}\widetilde G_n(t).
\]
Because \(\widetilde G_n(t)\xrightarrow{p}G_\ast(t)\), we have
\(\widetilde G_n(t)=O_p(1)\), and therefore
\[
\widehat G_n(t)=O_p(n^{-1}).
\]
\end{proof}

\begin{remark}[Explicit CLT variance of the scaled Greenwood estimator]
Using
\[
v_S(t)=S(t)^2h(t)
\]
and the explicit cross-covariance from the preceding remark, the delta-method
variance can also be written as
\[
\begin{aligned}
v_G(t)
&=
4S(t)^2h(t)^2v_S(t)
+S(t)^4v_H(t)
+4S(t)^3h(t)c_{SH}(t)\\
&=
S(t)^4\left[
v_H(t)
+2h(t)^3
-4h(t)\int_0^t\frac{\lambda(u)}{y(u)^2}\,du
\right].
\end{aligned}
\]
The expansion is nonnegative since \(v_G(t)=\operatorname{Var}\{\phi_{G,t}(Z)\}\).
This is the limiting-normal variance; no convergence of
\(n^3\operatorname{Var}\{\widehat G_n(t)\}\) to \(v_G(t)\) is asserted. No positivity condition on \(h(t)\) is
required for this product-map delta method; when \(h(t)=0\), the first-order
limit is degenerate, consistently with Section~\ref{subsec:consist}.
As an independent consistency check, the no-censoring specialization of this
formula is evaluated in Appendix~\ref{app:no-censoring}.
\end{remark}

\section{Discussion}

The results of this paper identify the natural $O_p(1)$ scale at which
Greenwood's variance estimator admits ordinary large-sample analysis.
Although Greenwood's estimator is usually written in the form
\[
\widehat G_n(t)
=
\widehat S_n(t)^2\,\widehat H_n(t),
\]
the Greenwood component \(\widehat H_n(t)\) itself is of order \(n^{-1}\).
Consequently,
\[
\widehat H_n(t)\xrightarrow{P}0,
\qquad
\widehat G_n(t)\xrightarrow{P}0,
\]
and the natural $O_p(1)$-scaled quantities are
\[
\widetilde H_n(t)=n\widehat H_n(t),
\qquad
\widetilde G_n(t)=n\widehat G_n(t).
\]

The main mathematical contribution of the paper is the pointwise
asymptotically linear representation, and the resulting central limit theorem,
for the scaled Greenwood component. Its deterministic limit is
\[
h(t)=\int_0^t \frac{\lambda(s)}{y(s)}\,ds,
\]
and
\[
\sqrt n\{\widetilde H_n(t)-h(t)\}
=
\frac1{\sqrt n}\sum_{i=1}^n\phi_{H,t}(Z_i)+o_p(1)
\xRightarrow{d}
\mathcal N\{0,v_H(t)\},
\]
with
\[
v_H(t)=\operatorname{Var}\{\phi_{H,t}(Z)\}.
\]
The scaling and consistency themselves are not the new first-order
contribution \cite[Theorem~4.2.2]{gill1980}; the analysis begins at the $\sqrt n$-scale around \(h(t)\).
The influence-function formulation shows why this refinement matters: the
first-order fluctuation contains both the martingale fluctuation of the event
counting process and the contribution of the empirical at-risk process. The
latter therefore enters the full first-order variance and cannot in general be
omitted.

Combining this representation with the standard first-order representation of
the Kaplan--Meier estimator \cite{reid1981,loSingh1986} gives the joint theory for
\((\widehat S_n(t),\widetilde H_n(t))\). In particular, their asymptotic
cross-covariance is generally nonzero and must be retained when the product map
is linearized. For the scaled Greenwood estimator,
\[
G_\ast(t):=S(t)^2h(t),
\qquad
\widetilde G_n(t)=\widehat S_n(t)^2\widetilde H_n(t),
\]
the resulting influence function is
\[
\phi_{G,t}(Z)
=
2S(t)h(t)\psi_{S,t}(Z)
+
S(t)^2\phi_{H,t}(Z),
\]
so that
\[
\sqrt n\{\widetilde G_n(t)-G_\ast(t)\}
=
\frac1{\sqrt n}\sum_{i=1}^n\phi_{G,t}(Z_i)+o_p(1)
\xRightarrow{d}
\mathcal N\{0,v_G(t)\}.
\]
Thus the natural object for consistency and first-order asymptotic normality is
not the unscaled Greenwood estimator \(\widehat G_n(t)\), but its scaled version
\(\widetilde G_n(t)\).

The case \(h(t)=0\) is included in the additive first-order theory. Then
\(\Lambda(t)=0\), the relevant influence functions vanish almost surely, and
\[
v_H(t)=v_G(t)=0.
\]
The corresponding central limit theorems therefore remain valid in degenerate
form. Positivity of \(h(t)\) is required only for transformations such as the
logarithm. When \(h(t)>0\), the same first-order theory yields the log-scale
limits for \(\widehat H_n(t)\); multiplying by the deterministic factor \(n\)
changes the centering by \(\log n\) but does not change the asymptotic log-scale
variance.

The assumptions used in the paper separate several distinct issues. Under the
i.i.d.\ sampling and independent-censoring structure in (A1), together with the
event-hazard regularity in (A3), the observed failure process has the
multiplicative-intensity representation used in the martingale expansion
\cite{aalen1978,fleming1991,andersen1993}.
Assumption (A2), combined with the uniform law of large numbers for the
empirical at-risk process, implies that the empirical risk set is bounded away
from zero uniformly on the time interval under consideration with probability
tending to one. Assumption (A4) is not needed for the first-order expansion of
the Nelson--Aalen plug-in term itself; it is used to simplify the exact
Greenwood--Nelson--Aalen remainder and to make that remainder negligible when
the expansion is transferred to \(\widetilde H_n(t)\).

The continuous-time no-ties assumption is used as a structural simplification
rather than as a claim that tied event times cannot occur in applications. It
is appropriate when the event-time model is continuous and ties are absent, or
when the no-ties model is regarded as an idealization of underlying continuous
event times whose observed ties arise only from recording precision. It is not,
however, without loss of generality for genuinely grouped or discrete
event-time data. The exact decomposition in Section~2.3 already shows how tied
failure multiplicities enter the Greenwood--Nelson--Aalen remainder. If the
multiplicities at observed failure times remain uniformly bounded, the same
type of denominator control suggests that this remainder should remain
negligible at both the consistency and central-limit scales. This observation
concerns only the Greenwood--Nelson--Aalen remainder; a full extension to
genuinely grouped or discrete event times would also require revisiting the
event-process representation and its first-order expansion. By contrast, if
some tied failure counts \(d_j\) are of order \(n\), as may occur with genuinely
grouped event times, the remainder term need not be negligible and the
continuous-time limit theory developed here no longer applies directly. Such
settings require a separate extension of the asymptotic expansion.

A further possible extension would be to develop uniform-in-time versions of
the consistency and central limit results, rather than the fixed-time results
considered here. Functional large-sample theory for the Kaplan--Meier process
provides an important precedent \cite{breslow1974,hallWellner1980,gill1983,majorRejto1988},
but extending the Greenwood-component expansion uniformly would still require
a functional central limit argument and corresponding tightness or
stochastic-equicontinuity control beyond the pointwise theory developed here.

\appendix
\section{Technical Proof of the Empirical Cross-Term Bound}
\label{app:cross-term-proof}

This appendix supplies the second-moment argument used in
Lemma~\ref{lem:empirical-cross-term}. Fix \(t\le\tau\) and suppose
Assumptions~\textup{(A1)}--\textup{(A3)} hold, as in that lemma. Recall that
\[
\nu(ds):=y(s)\lambda(s)\,ds.
\]
Put
\[
A_i(s):=Y_i(s)-y(s),
\qquad
\widetilde\nu_i(ds):=dN_i(s)-\nu(ds),
\]
and define
\[
K_{ij}:=K(Z_i,Z_j)
:=\int_0^t\frac{A_i(s)}{y(s)^3}\,\widetilde\nu_j(ds).
\]
Since
\[
y_n(s)-y(s)=\frac1n\sum_{i=1}^n A_i(s),
\qquad
\{\nu_n-\nu\}(ds)=\frac1n\sum_{j=1}^n\widetilde\nu_j(ds),
\]
we obtain
\[
C_n(t):=
\int_0^t\frac{y_n(s)-y(s)}{y(s)^3}\{\nu_n-\nu\}(ds)
=\frac1{n^2}\sum_{i=1}^n\sum_{j=1}^nK_{ij}.
\]
We split \(C_n=C_{n,\mathrm{diag}}+C_{n,\mathrm{off}}\) according as
\(i=j\) or \(i\ne j\).

Let \(m_t=\inf_{0\le s\le t}y(s)>0\). Moreover,
\[
\nu([0,t])
=\int_0^t y(s)\lambda(s)\,ds
\le \Lambda(t)<\infty.
\]
Since \(|A_i(s)|\le1\), \(N_j(t)\le1\), and the total variation
measure \(|\widetilde\nu_j|\) satisfies \(|\widetilde\nu_j|\le dN_j+\nu\),
\[
\begin{aligned}
|K_{ij}|
&\le
m_t^{-3}\int_0^t |A_i(s)|\,d|\widetilde\nu_j|(s)\\
&\le
m_t^{-3}\{N_j(t)+\nu([0,t])\}\\
&\le
m_t^{-3}\{1+\nu([0,t])\}
=:C_t<\infty
\end{aligned}
\]
for all \(i,j\). Consequently,
\[
|C_{n,\mathrm{diag}}|
\le \frac1{n^2}\sum_{i=1}^n|K_{ii}|
\le \frac{C_t}{n}.
\]

For \(i\ne j\), independence and centering give
\[
E(K_{ij}\mid Z_i)=0,
\qquad
E(K_{ij}\mid Z_j)=0.
\]
Indeed,
\[
E\{A_i(s)\}=0,
\qquad
E\Bigl\{\int_B dN_j(s)\Bigr\}=\nu(B),
\quad B\in\mathcal B([0,t]).
\]
The second identity is an equality of mean measures, since
\(P(X_j\in B,\Delta_j=1)=\int_B y(s)\lambda(s)\,ds=\nu(B)\).
Independence and the uniform total-variation bound justify conditional
Fubini, given \(Z_i\) (centered event measure) or \(Z_j\) (centered
risk indicator integrated against the fixed signed measure). Thus,
\[
E(K_{ij}\mid Z_i)
=
\int_0^t\frac{A_i(s)}{y(s)^3}
E\{\widetilde\nu_j(ds)\}
=0,
\]
and
\[
E(K_{ij}\mid Z_j)
=
\int_0^t\frac{E\{A_i(s)\}}{y(s)^3}
\widetilde\nu_j(ds)
=0.
\]
Thus the off-diagonal kernel is canonical in both arguments.

Now
\[
E\{C_{n,\mathrm{off}}^2\}
=\frac1{n^4}
\sum_{i\ne j}\sum_{k\ne l}E(K_{ij}K_{kl}).
\]
Call an index value among \(i,j,k,l\) a \emph{singleton index} if it occurs in
exactly one of the four positions. If such a singleton is present, condition on
all observations except the one carrying that singleton index. The factor that
depends on the remaining observation then has conditional expectation zero by
the canonical centering above, and therefore
\[
E(K_{ij}K_{kl})=0.
\]
Subject to \(i\ne j\) and \(k\ne l\), absence of singleton indices leaves only
\[
(i,j)=(k,l)
\qquad\text{or}\qquad
(i,j)=(l,k).
\]
There are at most \(2n(n-1)=O(n^2)\) such terms, and each satisfies
\[
|E(K_{ij}K_{kl})|
\le E|K_{ij}K_{kl}|
\le C_t^2.
\]
Thus
\[
E\{C_{n,\mathrm{off}}^2\}=O(n^{-2}).
\]
Markov's inequality applied to \(C_{n,\mathrm{off}}^2\) therefore gives
\[
C_{n,\mathrm{off}}=O_p(n^{-1}).
\]
Together with the deterministic diagonal bound this gives
\[
C_n(t)=O_p(n^{-1}),
\]
as asserted in Lemma~\ref{lem:empirical-cross-term}. Notice that this argument
uses only Assumptions~\textup{(A1)}--\textup{(A3)}; the no-ties
Assumption~\textup{(A4)} is not needed for the empirical cross-term bound.

\section{Expanded Variance of the Scaled Greenwood Component}
\label{app:expanded-H-clt}

Theorem~\ref{thm:CLT_tildeH} is stated in influence-function form in the main
text. The same central limit theorem can equivalently be written with the
variance expanded in integral form. This expanded representation separates the
martingale contribution from the additional terms induced by the randomness of
the empirical at-risk process.

\begin{proposition}[Expanded-variance formulation of Theorem~\ref{thm:CLT_tildeH}]
Under the assumptions of Theorem~\ref{thm:CLT_tildeH}, for each fixed
\(t\in[0,\tau]\),
\[
\sqrt n\{\widetilde H_n(t)-h(t)\}
\xRightarrow{d}
\mathcal N\{0,v_H(t)\},
\]
where
\[
h(t)=\int_0^t \frac{\lambda(s)}{y(s)}\,ds
\]
and
\begin{align}
v_H(t)
={}&
\int_0^t \frac{\lambda(s)}{y(s)^3}\,ds
\notag\\
&+
\int_0^t\int_0^t
\frac{\lambda(s)\lambda(r)}
     {y(s)^2y(r)^2}
\{y(s\vee r)-y(s)y(r)\}\,ds\,dr
\notag\\
&+
2\int_0^t
\frac{\lambda(r)}{y(r)}
\left\{
\int_0^r \frac{\lambda(s)}{y(s)^2}\,ds
\right\}
dr.
\label{eq:vH-expanded-app}
\end{align}
\end{proposition}

\begin{proof}
By Theorem~\ref{thm:CLT_tildeH},
\[
\sqrt n\{\widetilde H_n(t)-h(t)\}
\xRightarrow{d}
\mathcal N\{0,\operatorname{Var}(\phi_{H,t}(Z))\},
\]
where
\[
\phi_{H,t}(Z)
=
\int_0^t \frac{1}{y(s)^2}\,dM(s)
-
\int_0^t
\frac{\lambda(s)}{y(s)^2}
\{Y(s)-y(s)\}\,ds.
\]
It remains only to expand
\(\operatorname{Var}\{\phi_{H,t}(Z)\}\).

Set
\[
A_t
:=
\int_0^t \frac{1}{y(s)^2}\,dM(s),
\qquad
B_t
:=
\int_0^t
\frac{\lambda(s)}{y(s)^2}
\{Y(s)-y(s)\}\,ds.
\]
Then
\[
\phi_{H,t}(Z)=A_t-B_t.
\]
Let
\[
m_t:=\inf_{0\le u\le t}y(u)>0.
\]
Since \(\Lambda(t)<\infty\),
\[
\int_0^t\frac{\lambda(s)}{y(s)^2}\,ds
\le m_t^{-2}\Lambda(t)<\infty.
\]
Together with \(|Y(s)-y(s)|\le1\), this gives
\[
|B_t|\le m_t^{-2}\Lambda(t),
\]
so \(B_t\) is square integrable. The corresponding square integrability of
\(A_t\) follows from the martingale-isometry calculation below. These bounds
also justify the expectation--integration interchanges used in the variance
calculations.

First, \(d\langle M\rangle(s)=Y(s)\lambda(s)\,ds\). Hence, by the
martingale isometry and \(E\{Y(s)\}=y(s)\),
\[
\begin{aligned}
E(A_t^2)
&=
E\left[
\int_0^t \frac{1}{y(s)^4}\,d\langle M\rangle(s)
\right]\\
&=
\int_0^t \frac{\lambda(s)}{y(s)^3}\,ds.
\end{aligned}
\]
This is the martingale contribution.

Second, for any \(s,r\le t\),
\[
Y(s)Y(r)
=
\mathbf 1_{\{X\ge s\}}\mathbf 1_{\{X\ge r\}}
=
\mathbf 1_{\{X\ge s\vee r\}},
\]
and therefore
\[
E\{Y(s)Y(r)\}
=
y(s\vee r).
\]
Hence
\[
\operatorname{Cov}\{Y(s),Y(r)\}
=
y(s\vee r)-y(s)y(r),
\]
which gives
\[
\operatorname{Var}(B_t)
=
\int_0^t\int_0^t
\frac{\lambda(s)\lambda(r)}
     {y(s)^2y(r)^2}
\{y(s\vee r)-y(s)y(r)\}\,ds\,dr.
\]

Third, for fixed \(r\le t\), we verify the covariance identity needed below.
Under the continuous-time assumption, the probability of an observed failure
at the fixed endpoint \(r\) is zero, so the endpoint convention at \(r\) does
not affect the following decomposition. Split the martingale integral into its
part up to \(r\) and its subsequent increment:
\[
\int_0^t \frac{1}{y(s)^2}\,dM(s)
=
\int_0^r \frac{1}{y(s)^2}\,dM(s)
+
\int_r^t \frac{1}{y(s)^2}\,dM(s).
\]
Because \(Y(r)\) is \(\mathcal F_r\)-measurable, while conditional on
\(\mathcal F_r\) the subsequent martingale increment has conditional mean zero,
\[
E\left[
Y(r)\int_r^t \frac{1}{y(s)^2}\,dM(s)
\right]
=0.
\]
For \(s<r\), the single-event structure gives
\[
Y(r)\,dN(s)=0,
\qquad
Y(r)Y(s)=Y(r).
\]
Indeed, if a failure occurs at \(s<r\), then \(X=s<r\) and hence
\(Y(r)=0\); conversely, if \(Y(r)=1\), then \(X\ge r>s\), so also
\(Y(s)=1\). Thus, using
\(dM(s)=dN(s)-Y(s)\lambda(s)\,ds\),
\[
\begin{aligned}
E\left[
Y(r)\int_0^r \frac{1}{y(s)^2}\,dM(s)
\right]
&=
-\int_0^r
\frac{\lambda(s)}{y(s)^2}
E\{Y(r)Y(s)\}\,ds
\\
&=
-y(r)\int_0^r
\frac{\lambda(s)}{y(s)^2}\,ds.
\end{aligned}
\]
Therefore,
\[
E\left[
Y(r)\int_0^t \frac{1}{y(s)^2}\,dM(s)
\right]
=
-y(r)\int_0^r \frac{\lambda(s)}{y(s)^2}\,ds.
\]
Since \(E(A_t)=0\), it follows that
\[
\begin{aligned}
\operatorname{Cov}(A_t,B_t)
&=
\int_0^t
\frac{\lambda(r)}{y(r)^2}
E\{A_t(Y(r)-y(r))\}\,dr
\\
&=
-\int_0^t
\frac{\lambda(r)}{y(r)}
\left\{
\int_0^r \frac{\lambda(s)}{y(s)^2}\,ds
\right\}
dr.
\end{aligned}
\]
The negative sign is important: when this covariance is inserted into
\(\operatorname{Var}(A_t-B_t)\), the term \(-2\operatorname{Cov}(A_t,B_t)\)
produces the positive third term in \eqref{eq:vH-expanded-app}.
Combining
\[
\operatorname{Var}(A_t-B_t)
=
\operatorname{Var}(A_t)
+
\operatorname{Var}(B_t)
-
2\operatorname{Cov}(A_t,B_t)
\]
gives exactly \eqref{eq:vH-expanded-app}.

In particular, the term
\[
\int_0^t \frac{\lambda(s)}{y(s)^3}\,ds
\]
is the martingale contribution to the variance. In general, it is not the
full asymptotic variance: the risk-set contribution and its covariance
with the martingale term must also be retained. When \(h(t)=0\), all
three contributions vanish.

\end{proof}

\section{The No-Censoring Case as an Internal Consistency Check}
\label{app:no-censoring}

As an internal consistency check on the general first-order theory developed
above, we consider the special case of no censoring. In this setting the
risk-set process simplifies substantially, and several of the general
quantities admit closed-form expressions. This makes it possible to verify
independently the influence functions, the cross-covariance, and the
asymptotic variances derived for both the scaled Greenwood component and the
scaled Greenwood estimator.

Assume that the censoring time is almost surely infinite,
\[
P(C=\infty)=1.
\]
Then \(X=T\) almost surely and, under the maintained continuous-time model,
\[
y(t)=P(T\ge t)=P(T>t)=S(t).
\]
Consequently,
\[
h(t)
=
\int_0^t \frac{\lambda(s)}{S(s)}\,ds
=
\frac{1}{S(t)}-1,
\]
because \(S(0)=1\) and
\[
\frac{d}{ds}\{S(s)^{-1}\}
=
\frac{\lambda(s)}{S(s)}.
\]

The no-censoring case also gives an exact finite-sample identity. Let
\(k=N_n(t)\) be the number of failures observed by time \(t\). On the regular
event \(\mathcal D_n(t)=\{k<n\}\), the successive risk-set sizes at these
failures are \(n,n-1,\ldots,n-k+1\). Therefore,
\[
\begin{aligned}
\widehat H_n(t)
&=
\sum_{j=1}^k
\frac{1}{(n-j+1)(n-j)}\\
&=
\frac{1}{n-k}-\frac{1}{n},
\end{aligned}
\]
while
\[
\widehat S_n(t)=\frac{n-k}{n}.
\]
Hence, exactly on \(\mathcal D_n(t)\),
\begin{equation}\label{eq:no-censoring-H-identity}
\boxed{
\widetilde H_n(t)
=
n\widehat H_n(t)
=
\frac{1}{\widehat S_n(t)}-1.
}
\end{equation}
Since \(S(t)>0\) under the standing assumptions,
\[
P\{\mathcal D_n(t)^c\}
=
P\{N_n(t)=n\}
=
\{1-S(t)\}^n
\longrightarrow0,
\]
so the arbitrary global extension of \(\widehat H_n(t)\) outside
\(\mathcal D_n(t)\) does not affect the following first-order conclusions.

Without censoring, \(\widehat S_n(t)\) is simply the empirical survival
probability,
\[
\widehat S_n(t)
=
\frac1n\sum_{i=1}^n Y_i(t),
\]
and therefore
\[
\sqrt n\{\widehat S_n(t)-S(t)\}
\xRightarrow{d}
\mathcal N\bigl(0,S(t)\{1-S(t)\}\bigr),
\qquad
\psi_{S,t}(Z)=Y(t)-S(t).
\]
Applying the one-dimensional delta method to
\[
f(x)=x^{-1}-1,
\qquad
f'\{S(t)\}=-S(t)^{-2},
\]
and using the exact identity \eqref{eq:no-censoring-H-identity} gives the
special-case influence function
\begin{equation}\label{eq:no-censoring-phiH}
\boxed{
\phi_{H,t}(Z)
=
-\frac{Y(t)-S(t)}{S(t)^2}
=
\frac{S(t)-Y(t)}{S(t)^2}.
}
\end{equation}
Thus
\begin{equation}\label{eq:no-censoring-vH}
\boxed{
v_H(t)
=
\frac{1-S(t)}{S(t)^3}.
}
\end{equation}
This agrees with the no-censoring specialization of the expanded variance
formula \eqref{eq:vH-expanded-app} and therefore provides a direct check of
the influence-function result for \(\widetilde H_n(t)\).

By contrast, retaining only the martingale part of the general influence
function would give
\[
\int_0^t \frac{\lambda(s)}{S(s)^3}\,ds
=
\frac13\{S(t)^{-3}-1\},
\]
because
\[
\frac{d}{ds}\{S(s)^{-3}\}
=
3\lambda(s)S(s)^{-3}.
\]
This is generally different from \eqref{eq:no-censoring-vH}. Thus the
no-censoring case already shows explicitly that the empirical risk-set
contribution cannot in general be omitted from the first-order variance.

The same specialization also checks the cross-covariance in
Section~\ref{subsec:joint-AN}. From \eqref{eq:no-censoring-phiH},
\[
\begin{aligned}
c_{SH}(t)
&=
\operatorname{Cov}\{\psi_{S,t}(Z),\phi_{H,t}(Z)\}\\
&=
-\frac{1}{S(t)^2}
\operatorname{Var}\{Y(t)-S(t)\}\\
&=
-\frac{1-S(t)}{S(t)}.
\end{aligned}
\]
On the other hand,
\[
J(t)
:=
\int_0^t\frac{\lambda(s)}{S(s)^2}\,ds
=
\frac12\{S(t)^{-2}-1\},
\]
since
\[
\frac{d}{ds}\{S(s)^{-2}\}
=
2\lambda(s)S(s)^{-2}.
\]
Substituting this and \(h(t)=S(t)^{-1}-1\) into the general formula
\[
c_{SH}(t)
=
-S(t)\left\{J(t)+\frac12h(t)^2\right\}
\]
gives exactly the same value,
\[
c_{SH}(t)=-\frac{1-S(t)}{S(t)}.
\]

Finally, \eqref{eq:no-censoring-H-identity} yields another exact
finite-sample identity on \(\mathcal D_n(t)\):
\begin{equation}
\boxed{
\widetilde G_n(t)
=
\widehat S_n(t)^2\widetilde H_n(t)
=
\widehat S_n(t)\{1-\widehat S_n(t)\}.
}
\end{equation}
With
\[
q(x)=x(1-x),
\qquad
q'\{S(t)\}=1-2S(t),
\]
the delta method therefore gives
\begin{equation}
\boxed{
\phi_{G,t}(Z)
=
\{1-2S(t)\}\{Y(t)-S(t)\},
}
\end{equation}
and hence
\begin{equation}\label{eq:no-censoring-vG}
\boxed{
v_G(t)
=
S(t)\{1-S(t)\}\{1-2S(t)\}^2.
}
\end{equation}
The general influence-function formula from Corollary~\ref{cor:scaled-greenwood}
reduces to the same expression, because
\[
2S(t)h(t)=2\{1-S(t)\},
\qquad
S(t)^2\phi_{H,t}(Z)=-\psi_{S,t}(Z),
\]
so that
\[
\phi_{G,t}(Z)
=
2S(t)h(t)\psi_{S,t}(Z)
+S(t)^2\phi_{H,t}(Z)
=
\{1-2S(t)\}\psi_{S,t}(Z).
\]

As a further algebraic check, substitute
\[
h(t)=\frac{1-S(t)}{S(t)},
\qquad
J(t)=\frac12\{S(t)^{-2}-1\},
\qquad
v_H(t)=\frac{1-S(t)}{S(t)^3}
\]
into the expanded variance formula from Section~\ref{subsec:delta-method-consequences},
\[
v_G(t)
=
S(t)^4\{v_H(t)+2h(t)^3-4h(t)J(t)\}.
\]
Then
\[
\begin{aligned}
v_H(t)+2h(t)^3-4h(t)J(t)
&=
\frac{1-S(t)}{S(t)^3}
\left[
1+2\{1-S(t)\}^2-2\{1-S(t)^2\}
\right]\\
&=
\frac{1-S(t)}{S(t)^3}\{1-2S(t)\}^2,
\end{aligned}
\]
which again yields \eqref{eq:no-censoring-vG}.

When \(h(t)=0\), the no-censoring identity
\(h(t)=S(t)^{-1}-1\) implies \(S(t)=1\). In that case
\(\phi_{H,t}=\phi_{G,t}=0\) almost surely and
\(v_H(t)=v_G(t)=0\), consistently with the degenerate cases discussed in the
main text.

Thus, in the no-censoring case, the general influence-function
representations, the cross-covariance formula, and the asymptotic variances
for both \(\widetilde H_n(t)\) and \(\widetilde G_n(t)\) reduce to the
corresponding expressions obtained directly from the elementary identities
above. This provides an independent internal consistency check of the
first-order theory developed in Sections~\ref{sec:first-order-H} and~4.

\newpage
\section*{Statements and Declarations}

\noindent\textbf{Competing interests.} The author declares no competing interests.

\noindent\textbf{Funding.} No specific funding was received for this work.

\noindent\textbf{Data availability.} No data were generated or analysed in this study.

\noindent\textbf{Use of generative AI.} The research question, the underlying mathematical concept, several key proof ideas, and substantial parts of the mathematical derivations and proofs were developed by the author. Generative artificial intelligence (ChatGPT, OpenAI) was used to assist in the further development and verification of mathematical arguments, the examination of proof details and gaps, and the preparation and systematic review of the manuscript. The author critically evaluated the AI-assisted material and assumes full responsibility for the scientific content and conclusions of this work.

\end{document}